\documentclass[11pt,a4paper]{amsart}
\usepackage{amssymb}
\usepackage{epsfig}
\begin{document}
\newtheorem{thrm}{Theorem}%[chapter]
\newtheorem{thmf}{Th{\'e}or{\`e}me}%[chapter]
\newtheorem{thmm}{Theorem}%[chapter]
\def\thethmm{\ref{twist}$'$}
\newenvironment{thint}[1]{{\flushleft\sc{Th{\'e}or{\`e}me}}
      {#1}. \it}{\medskip} 
\newenvironment{thrm.}{{\flushleft\bf{Theorem}}. \it}{\medskip} 
\newenvironment{thrme}[1]{{\flushleft\sc{Theorem}}
      {#1}. \it}{\medskip} 
\newenvironment{propint}[1]{{\flushleft\sc{Proposition}}
      {#1}. \it}{\medskip} %\def\thethrme{\ref{cls}}
\newenvironment{corint}[1]{{\flushleft\sc{Corrolaire}}
      {#1}. \it}{\medskip} %\def\thethrme{\ref{cls}}
\newtheorem{corr}{Corollary}%[chapter]
\newtheorem{corf}{Corollaire}%[chapter]
\newtheorem{prop}{Proposition}%[chapter]
\newtheorem{defi}{Definition}%[chapter]
\newtheorem{deff}{D{\'e}finition}%[chapter]
\newtheorem{lem}{Lemma}%[chapter]
\newtheorem{lemf}{Lemme}%[chapter]
\def\fy{\varphi}
\def\ul{\underline}
\def\obsf{{\flushleft\bf Remarque. }}
\def\obs{{\flushleft\bf Remark. }}
\def\R{\mathbb{R}}
\def\C{\mathbb{C}}
\def\Z{\mathbb{Z}}
\def\N{\mathbb{N}}
\def\P{\mathbb{P}}
\def\Q{\mathbb{Q}}
\def\p{\pi_1(X)}
\def\Re{\mathrm{Re}}
\def\Im{\mathrm{Im}}
\def\H{\mathbf{H}}
\def\Hs{{Nil}^3}
\def\L{L_\alpha}
\def\h{\frac{1}{2}}
\def\M{\P(E)\times\P(E)^*\smallsetminus\mathcal{F}}
\def\m{\cp2\times{\cp2}^*\smallsetminus\mathcal{F}}
\def\g{\mathfrak{g}}
\def\l{\mathcal{L}}
\def\k{\mathfrak{h}}
\def\s3{\mathfrak{s}^3}
\def\nb{\nabla}
\def\Re{\mathrm{Re}}
\def\Im{\mathrm{Im}}
\newcommand\cp[1]{\mathbb{CP}^{#1}}
\renewcommand\o[1]{\mathcal{O}({#1})}
\renewcommand\d[1]{\partial_{#1}}
\def\sq{\square_X}
\def\gp{\dot\gamma}
\def\j{\mathcal{J}}
\def\iif{\mbox{\bf\em{I}}}

\title{Normal $CR$ structures on compact 3-manifolds}

\author{Florin Alexandru Belgun}

\begin{abstract} We study normal $CR$ compact manifolds in dimension 3. For a
  choice of a $CR$ Reeb vector field, we associate a Sasakian metric
  on them, and we classify those metrics. As a consequence, the underlying
  manifolds are topologically finite quotiens of $S^3$ or of a circle
  bundle over a Riemann surface of positive genus. In the latter case,
  we prove that their  $CR$ automorphisms group is a finite extension
  of $S^1$, and we classify the normal $CR$ structures on these manifolds.
\end{abstract}

\maketitle

\section{Introduction}\label{s.k1} %From the differential-geometric point of view,
%one of the nicest geometric structure one can expect to find on an
%even-dimensional manifold is a complex structure. There is no perfect
%analog of it on odd-dimensional manifolds, but the concept of $CR$
%structure arises naturally from the geometry of real hypersurfaces of
%complex spaces, e.g. the boundaries of pseudo-convex domains in
%$\C^n$. Perhaps the first imperfection comes from the non-ellipticity of the
%Cauchy-Riemann operator on a $CR$ manifold, such that analyticity no
%longer follows from formal integrability (see next Section).

Analogs of complex manifolds in odd dimensions, pseudo-conformal $CR$
manifolds are particular contact manifolds, with a complex structure
on the corresponding distribution of hyperplanes, satisfying an
integrability condition (see Section \ref{s.k2}). Contrary to complex
geometry, $CR$ geometry is locally determined by a finite system of
local invariants (like in the cases of conformal or projective
structures), \cite{tk1}, \cite{kob}, \cite{stern}. Therefore the space
of locally non-isomorphic $CR$ structures is a space with infinitely
many parameters.

In this paper, we focus our attention on {\it normal} $CR$ manifolds,
which admit global {\it Reeb} vector fields preserving the $CR$
structure, in particular their $CR$ automorphisms group has dimension
at least 1. Our main result is that, for a compact normal $CR$
3-manifold, which is topologically not a quotient of $S^3$, this $CR$
automorphisms group is a finite extension of a circle, thus the Reeb
vector field is unique up to a constant (Section 4, Theorem
\ref{t2}). This, together with the classification of {\it Sasakian}
compact 3-manifolds (see Section 3),  allows us to obtain the
classification of normal $CR$ structures on these manifolds (Section
4, Corollary \ref{cor1}).

The question of classifying compact $CR$ manifolds has first been
solved in situations with a high order of local symmetry: the
classification of {\it flat} compact $CR$ manifolds, where the local
$CR$ automorphism group is $PSU(n+1,1)$ (if the manifold has dimension
$2n+1$), is due to E. Cartan \cite{cart} and to D. Burns and
S. Shnider \cite{bsh}; in dimension 3,  homogeneous, simply-connected,
$CR$ manifolds are either flat or (3-dimensional) Lie groups, and have
been classified by E. Cartan \cite{cart} (see also \cite{ens}). In
this case, there is no intermediate symmetry because E. Cartan has
showed that a homogeneous $CR$ manifold whose $CR$ automorphism group
has dimension greater than 3 is automatically flat.
%\iVs\isk\ipscf\icrn\iS

In dimension 3, the normal $CR$ structures are always deformations of
a flat one (Theorem \ref{sasaki}, see also \cite{F}), and the key
point is that, for a $CR$ Reeb vector field $T$, they admit compatible
{\it Sasakian metrics}, for which $T$ is Killing (see Section
\ref{s.k2} for details); these metrics are closely related to locally
conformally K{\"a}hler metrics with parallel Lee form, natural analogs
of K{\"a}hler structures on non-symplectic complex manifolds
\cite{vais}.
% In this paper, we consider, in dimension 3, {\it normal} $CR$
% manifolds (see Section \ref{s.k2}), which admit at least one
% infinitesimal $CR$ 
% automorphism (a {\it  CR Reeb vector field}), transversal to the complex
% (hyper)plane. We show that, on each compact 
% 3-manifold, such structures are deformations of a flat one (Theorem
% \ref{sasaki}, see also \cite{F}), so the present study (of a far less
% symmetric situation, as we shall see in Section \ref{s.k4}) follows
% naturally after the previously cited works. A choice of a $CR$ Reeb
% vector field $T$ yields a {\it Sasakian metric} on the manifold, for which
% $T$ is a Killing vector field (see
% Section \ref{s.k2} for details); these metrics are closely related to locally
% conformally K{\"a}hler metrics with parallel Lee form, natural analogs
% of K{\"a}hler structures on non-symplectic complex manifolds
% \cite{vais}.

Topologically, every compact normal $CR$ (or Sasakian) 3-manifold is a
Seifert fibration (Pro\-position \ref{p1sk}, see also \cite{kamish},
\cite{gei} and \cite{F}), but it turns out that the Sasakian
structures themselves can be explicitly described on these manifolds:
Theorem \ref{sasaki}, Proposition \ref{p1sk} (these are extended
versions of the results announced in \cite{F}); In particular, if a
compact Sasakian 3-manifold is not covered by $S^3$ (or {\it
non-spherical}), its Sasakian structure is always {\it regular},
i.e. it is a (finite quotient of a) circle bundle over a Riemann
surface and its metric arises from a Kaluza-Klein construction
(Corollary \ref{KK}), which leads to an elementary description of any
Sasakian metric on these  manifolds, either directly, or as
deformations of a $CR$ flat one.

Although the classification of Sasakian structures on compact
3-manifolds is complete in all cases (the explicit description in the
case of finite quotiens of $S^3$ is more elaborate, but still
possible, see Section 3), the question of classifying the associated
$CR$ structures (which are normal) is more subtle. It is solved in
Section 4 for non-spherical manifolds, and is still open for $S^3$.
\medskip

{\sc Acknowledgements.} The author is grateful to P. Gauduchon for his
constant support during the last few years.

\section{$CR$ geometry on 3-manifolds and Sasakian structures}\label{s.k2}

Let $M$ be a $2n+1$-dimensional manifold. In all generality, a $CR$
structure on $M$ is a field of complex structures on a field of
hyperplanes of $M$. But, as this concept is inspired by the structure
of a real hypersurface in a complex manifold, one generally searches
for $CR$ structures satisfying the conditions described below: 

  Let $Q$ be a field of hyperplanes in a $2n+1$-dimensional real
  manifold $M$, then the {\it Levi form} $\mathbf{L}^Q:Q\times
  Q\rightarrow TM/Q$ is defined by
%\index{Levi, forme de}
  $\mathbf{L}(X,Y):=[X,Y]\mod{Q}$. Let $J\in \mathrm{End}(Q)$,
  $J^2=-\mathbf{I}$ be an almost-complex structure on $Q$. We say that
  the Levi form is non-degenerate iff $\mathbf{L}^Q(X,\cdot)\in
  \mathrm{Hom}(Q,TM/Q)$ is non-zero for any non-zero $X\in Q$, and
  that it is $J$-invariant iff
  $\mathbf{L}^Q(J\cdot,J\cdot)=\mathbf{L}^Q(\cdot,\cdot)$. If the
  Levi form is $J$-invariant, the {\it
  Nijenhuis} {\it tensor} of $J$ is defined as a linear map
  $N:\Lambda^2Q\rightarrow Q$, by:
%\index{Nijenhuis, tenseur de!CR@$CR$}
$$4N(X,Y)=[JX,JY]-J[JX,Y]^Q-J[X,JY]^Q-[X,Y],$$
where $Z^Q$ denotes the component in $Q$ of the vector $Z$ --- by means
of a non-canonic linear projection ---; the Nijenhuis tensor is
independent of this projection.\smallskip
 
{\noindent\bf Convention } We call a {\it tensor} any multi-linear
object defined on subspaces or/and quotients of $TM$ (namely $Q$ and
$TM/Q$). For example, $\mathbf{L}^Q:Q\times
  Q\rightarrow TM/Q$ is a tensor.
\smallskip

\begin{defi}
Let $M$ be a odd-dimensional connected real manifold. A {\rm CR} structure on
$M$ is a field of hyperplanes $Q$, with an almost complex structure
$J\in\mathrm{End}(Q)$, such that the Levi form $\mathbf{L}^Q$ is
$J$-invariant.
The $CR$ structure $J$ on $M$ is called {\rm formally integrable} if the
Nijenhuis tensor vanishes identically; it is called pseudo-conformal
if the Levi form is non-degenerate, and pseudo-convex if
$\mathbf{L}^Q(J\cdot,\cdot)$ is a positive definite (symmetric)
2-form on $Q$.
\end{defi}
%\icr\icrfi\iSp

Some authors consider only (formally) integrable $CR$ structures; this
is because only these can arise as real hypersurfaces in a complex
manifold (in which case we call them {\it integrable}). Note that, in
the
$C^\infty$ case, the vanishing of the Nijenhuis tensor does not
necessarily imply integrability (the
analog of the Newlander-Nirenberg theorem holds only for analytic $CR$
manifolds)
\cite{nir1}, \cite{nir2}.
%\ipscf\icrfi\ipscv

We consider only pseudo-conformal $CR$ structures: 
pseudo-conformal geometry is, like projective or conformal
geometry, a {\it semisimple G-structure} (or a {\it parabolic
  geometry}) \cite{stern}, \cite{tk1}, in
particular they admit a unique {\it Cartan connection}, whose
curvature (see below) locally characterizes the geometry of a
pseudo-conformal 
manifold; another consequence of this is that the group of
diffeomorphisms of $M$ preserving a given pseudo-conformal structure
is a Lie group \cite{kob}.
%\iG\iC\iR\iTW

If $(M,Q,J)$ is a pseudo-conformal manifold, then $Q$ is a {\it
  contact structure} 
%\istc\ifc
on $M$, i.e. $Q$ is (locally) the kernel of a 
  1-form $\eta:TM\rightarrow L$ with values in a real line bundle,
  such that $\eta\wedge d\eta$ does not vanish. We consider only
  orientable manifolds $M$, such that $L$ is also orientable (hence
  topologically trivial), and then $\eta\wedge d\eta$ is a volume form
  on $M$. 
%\ifc
Obviously, to each vector field $X\in TM\smallsetminus Q$, we
  uniquely associate a contact form $\eta$ such that $\eta(X)\equiv
  1$. Conversely, to each contact form $\eta$ we can uniquely
  associate a {\it Reeb} vector field $T$ such that $\eta(T)=1$ and
  $d\eta(T,\cdot)=0$ everywhere (in particular, the Lie derivative of
  $\eta$ along $T$ vanishes: $\mathcal{L}_T\eta=0$).

The hyperplane $Q$ admits a natural
conformal- (pseudo-) Hermitian\linebreak
structure, represented by the (pseudo-)Hermitian symmetric forms
$h:=-\frac{1}{2}d\eta(J\cdot,\cdot)$ on $Q$, for any contact form
$\eta$ (the factor $\frac{1}{2}$ is useful when we consider Sasakian
metrics --- see below).  

The choice of a contact form $\eta$ (and of its Reeb vector field $T$)
yields the {\it Tanaka-Webster connection} $\nb$, defined as follows
\cite{tk2}, \cite{web}:
\begin{enumerate}
\item $T$, $Q$, and $J$ are parallel with respect to $\nb$;
\item If $\tau(X,Y):=\nb_XY-\nb_YX-[X,Y]$ is the torsion of $\nb$,
  then $\tau(X,Y)=d\eta(X,Y)T,\,\forall X,Y\in Q$;
\item If ${\tilde\tau}(X):=\tau(T,X), \, X\in Q$, then ${\tilde\tau}:Q\rightarrow Q$ is
  $J$-anti-invariant (${\tilde\tau}(JX)=-J{\tilde\tau}(X)$).
\end{enumerate}
We remark that this connection preserves the $CR$ structure, and it
has minimal torsion. ${\tilde\tau}$ cannot vanish unless $T$ is a $CR$ Reeb
vector field, i.e. the diffeomorphism group generated by $T$ on $M$
(that already preserves $Q$, as $T$ is a Reeb vector field) preserves
$J$, i.e. $\mathcal{L}_TJ=0$.

\begin{defi}
  The pseudo-conformal structure of $M$ is called {\rm normal} iff it
  admits a {\rm CR} Reeb vector field.
\end{defi}

In particular, the dimension of the Lie group of pseudo-conformal
automorphisms of a normal $CR$ manifold is at least 1.

If $T'$ is another Reeb vector field on $M$, then the corresponding
contact form $\eta'$ is equal to $f^{-1}\eta$, for a positive function
$f:M\rightarrow\R$, and 
$$T'=fT+X_f,$$ 
where $X_f\in Q$ is determined by
the fact that $\mathcal{L}_{T'}\eta'=0$, thus
\begin{equation}
  \label{eq:xf}
  d\eta(X_f,\cdot)=-df(X)\Leftrightarrow X_f=-\frac{1}{2}(df|_Q\circ
  J)^\sharp =\frac{1}{2}J(df|_Q)^\sharp, 
\end{equation}
where the ``rising of indices'' $df^\sharp$ is made with respect to
the Hermitian metric $h$ on $Q$. Then $\mathcal{L}_{T'}J=0$ iff
$[T'JX]-J[T',X]=0,\,\forall X\in Q$, thus
$$-\nb_{JX}T'-d\eta(X_f,JX)T+J(\nb_XT'+d\eta(X_f,X)T)=0,$$
where $\nb$ is the Tanaka-Webster connection corresponding to
$\eta$. We get then
$$-J\nb_{JX}(df|_Q)^\sharp-\nb_X(df|_Q)^\sharp=0,$$
which leads to (see also \cite{go}):
\begin{prop}
  \label{hess} If $f$ is a function relating 2 $CR$ Reeb vector fields
  $T$ and $T'=fT+X_f$, $X_f\in Q$, then the {\rm Hessian} of $f$,
  restricted to $Q$ (defined by
  $\mathrm{Hess}^Qf(X,Y):=X.Y.f-\nb_XY.f$), is $J$-invariant (with
  respect to the  Tanaka-Webster connection $\nb$ on $M$ induced by
  the contact form $\eta$): 
$$\mathrm{Hess}^Qf(X,Y)=\mathrm{Hess}^Qf(JX,JY),\,\forall X,Y\in Q.$$
In particular, if $\dim M=3$, the above condition means that
$\mathrm{Hess}^Qf$ is a scalar multiple of $h$.
\end{prop}
\subsection{Sasakian geometry}\label{ss.k21} A {\it Sasakian
structure} on an 
odd-dimensional manifold $M$ is a Riemannian metric $g$ on $M$ and a
unitary Killing vector field $T$ such that $\nb_TT=0$ and $\nb_\bullet
T:Q\rightarrow Q$ (where $Q:=T^\perp$) is an almost complex structure
$J$ (compatible with the metric as it is an anti-symmetric
endomorphism). It is easy to see that $T$ preserves $\nb T$, as it is
Killing, so a Sasakian structure on $M$ is a special case of a normal
(pseudo-convex) $CR$ structure on $M$ (we remark that $T$, followed by
a basis of $Q$, oriented by $J$, yield to an orientation of $M$).

Actually, if $\dim M=3$, then any $CR$ Reeb vector field $T$ of a
normal $CR$ structure on $M$ yields a unique Sasakian structure: 
\begin{prop}
  If $T$ is a $CR$ Reeb vector field on the 3-manifold $M$, then there
  is a unique Sasakian metric $g$ on $M$ for which $T$ or
  $-T$ is the corresponding Killing vector field.
\end{prop}
\begin{proof} As
$T$ preserves $Q$ and $J$, it (and all its multiples by a constant)
preserves the Riemannian metric defined by
\begin{equation}
  \label{eq:rbsk}
  g:=\eta^2-\frac{1}{2}d\eta(J\cdot,\cdot),
\end{equation}
(we replace, if necessary, $T$ with $-T$, $\eta$ with $-\eta$, such
that $g$ is positive definite) with respect to which $T$ is Killing
and $\nb_TT=0$. $\nb T$ is then identified to a $J$-invariant,
anti-symmetric endomorphism of $Q$, thus equal to $fJ$, for a function
$f$. But 
$$d\eta(X,Y)=(\nb_X\eta)(Y)-(\nb_Y\eta)(X)=-2g(fJX,Y),\;\forall X,Y\in
Q,$$
as $\eta=g(T,\cdot)$, thus $f\equiv 1$
\end{proof}

\obs For a general, pseudo-conformal manifold $M^{2n+1}$, let $(p,q),
\, p\leq q, \, p+q=n,$ be the signature of the Hermitian form
$h$. Then the {\it flat 
  model} of the pseudo-conformal geometry of signature $(p,q)$ is the
homogeneous space $PSU(p+1,q+1)/H^{p,q}$, where $H^{p,q}$ is the
isotropy subgroup 
of the point $[1:0:\ldots:0]\in\cp{n+1}$, for the action of
$PSU(p+1,q+1)$ on the real hypersurface $M^{p,q}$ of $\cp{n+1}$
defined by the equation
%\ipscf\isk\icrp\iC\iCc
$$x_1\bar x_1+\ldots +x_p\bar x_p-x_{p+1}\bar x_{p+1}
-\ldots-x_{p+q}\bar x_{p+q}+x_0\bar x_{n+1}=0.$$
It turns out that it exists a canonical $H^{p,q}$-bundle $P$ over any
pseudo-conformal manifold $M^{2n+1}$ (whose Hermitian structure on $Q$
has signature $(p,q)$), and a canonical {\it Cartan connection}
$\omega:TP\rightarrow \mathfrak{psu}(p+1,q+1)$ (where the latter is the
Lie algebra of the above mentioned group) \cite{tk1}, \cite{chm}. Its
{\it curvature} 
measures the obstruction to the construction of a local diffeomorphism
$\Xi:P\rightarrow PSU(p+1,q+1)$, for which $\omega$ would be the
differential (in particular, $\Xi$ would induce a group structure on
the universal covering of $P$, locally isomorphic to $PSU(p+1,q+1)$),
and it locally determines the pseudo-conformal structure, see Tanaka
\cite{tk1}; see also \cite{tk2}, \cite{chm}; see \cite{kob} for a
general theory of Cartan connections, and \cite{stern} for a general
theory of simple graded Lie algebras and $G$-structures.
\smallskip

The curvature of the Cartan connection is identified, if $n>1$, to the
pseudo-conformal tensor of Chern and Moser \cite{chm}, which is a
component of the curvature of any Tanaka-Webster connection
\cite{tk1}, \cite{web}. It is the equivalent of the Weyl tensor in
conformal geometry \cite{tk1}. In dimension 3, i.e. $n=1$, this tensor
vanishes identically and the curvature of the Cartan connection is
determined by another tensor (see the convention above), called the
{\it Tanaka curvature}, see \cite{tk1}, page 187, Theorem 12.3.

In the case when the Tanaka-Webster connection $\nb$ corresponds to a
{\it positive}
Reeb vector field $T$ (i.e. $h(X,X)=-\frac{1}{2}d\eta(JX,X)>0,\,\forall X\in Q,\,
X\neq 0$), we denote by $k$ the {\it sectional curvature} of
the plane $Q$: $k:=h(R(X,JX)X,JX)$, for $X\in Q$, $h(X,X)=1$. The
Tanaka curvature is then defined as the tensor $\Phi:S^2Q\rightarrow
TM/Q$ that satisfies:
%\iTWt\iTW\iT\isect
\begin{enumerate}
\item $\Phi$ is trace-free, i.e. $\Phi(X,X)+\Phi(JX,JX)=0$;
\item $\Phi(X,X)(T)=-kh({\tilde\tau}(X),X)+2h((\nb_T{\tilde\tau})(X),JX)-$\\
  $-\h(\nb_X\nb_{JX}-\nb_{\nb_XJX}+\nb_{JX}\nb_X-\nb_{\nb_{JX}X})k-$\\
$4h((\Delta{\tilde\tau})(X),X),$ \end{enumerate}
for any $X\in Q$; the {\it Laplacian} $\Delta$ is
  defined as
$$\Delta
\sigma:=(\nb_Y\nb_Y-\nb_{\nb_YY}+\nb_{JY}\nb_{JY}-\nb_{\nb_{JY}(JY)})\sigma,$$
for a unitary $Y\in Q$, and a tensor $\sigma\in
\mathrm{End}(Q)$. $\Phi$ is independent of the Reeb field $T$, and of
the associated connection $\nb$ \cite{tk1}. 
\obs The Tanaka curvature is the analog of the Cotton-York tensor from
conformal geometry~; however, the terms contained in the expression of
$\Phi$ change by terms involving up to 4th order derivatives of $f$,
for a change of the Tanaka-Webster connection determined by
$\eta'=f^{-1}\eta$. The invariance, in dimension 3, of the Tanaka
curvature is, thus, a highly non-trivial fact.
\smallskip

We denote by $\square_{X}$ the
second order differential operator on functions \linebreak $f:M\rightarrow\R$,
acting as 
\begin{equation}
  \label{eq:sq}
  \sq f:= X.JX.f+JX.X.f-\nb_XJX.f-\nb_{JX}X.f,
\end{equation}
and it is obvious that $\sq$ depends quadratically on $X\in Q$. Then,
if $T$ is a positive $CR$ Reeb vector field, then the above expression
for $\Phi$ becomes a lot simpler:
\begin{equation}
  \label{eq:fisk}
  \Phi(X,X)(T)=-\h\sq k,\; \forall X\in Q.
\end{equation}
\begin{prop}
  \label{k} The sectional curvature of the plane $Q$ in $M$, with
  respect to the Sasakian structure induced by $T$, is equal to $-2k-3$. 
\end{prop} 
\begin{proof} We denote by $\nb$, resp. $\nb^0$,  the Tanaka-Webster
connection, resp. the Levi-Civita (Sasakian) connection associated to the
positive $CR$ Reeb vector field $T$. We have 
\begin{eqnarray*}
  \nb_TT-\nb^0_TT&=&0;\\
\nb_XT-\nb^0_XT&=&-JX,\;\forall X\in Q;\\
\nb_TX-\nb^0_TX&=&-JX,\;\forall X\in Q;\\
\nb_XY-\nb^0_XY&=&\h Td\eta(X,Y)=Tg(JX,Y),\,\forall X,Y\in
Q.
\end{eqnarray*}
The claimed result now follows from a straightforward computation.
\end{proof}
On the other hand, if we define the operator $\sq^0$ as in
(\ref{eq:sq}), by replacing $\nabla$ with $\nabla^0$, we have
\begin{equation}
  \label{eq:sqsk}
  \sq f=\sq^0 f,\,\forall f:M\rightarrow\R,\, \forall X\in Q.
\end{equation}
\subsection{Regular Sasakian structures}\label{ss.k22} A Sasakian
structure on a 
compact 3-mani\-fold is called {\it regular} if the Reeb vector field
$T$ is induced by a free circle action on $M$. In
that case, there is a $S^1$ fibration
$M\stackrel{\pi}{\rightarrow}\Sigma$ of $M$ over a Riemann surface
$\Sigma$, $T$ is tangent to the fibers, and $Q$ is a connection in the
principal bundle $M\rightarrow \Sigma$, of connection form
$i\eta:TM\rightarrow i\R$ (where $i\R$ is the Lie algebra of
$S^1\simeq U(1)$), and of curvature form $id\eta$. On the other hand,
$T$ is a Killing vector field, thus the metric $g$ on $M$ induces a
Riemannian metric $g^0$ on $\Sigma$, compatible with the induced
complex structure $J$; the K{\"a}hler form on $\Sigma$ is then
$\omega=g^0(J\cdot,\cdot)$, thus $\pi^*\omega=\h d\eta$. Given such a
metric on $\Sigma$ and a connection, the metric constructed as above
(the horizontal space of the connection is defined to be orthogonal to
the vertical --- this procedure is usually called a {\it Kaluza-Klein
  construction}) on $M$ is Sasakian. 
%\iskr\iTW\icrnns\\icrp

\obs It turns out then that the Chern class of the $S^1$-bundle
$M\rightarrow\Sigma$ is always positive: this is because we chose $T$
to be {\it positive} (see \cite{F}, Section 3, for a detailed
explanation); in particular, we obtain a {\it positive} Chern class
for the Hopf fibration $S^3\rightarrow S^2$, apparently contradictory
to the {\it negative} Chern class of the tautological bundle
$\mathcal{O}(-1)$ on $\cp1$ ($\C^2\smallsetminus\{0\}\rightarrow\cp1$);
this is because the canonical metric on $S^3$ is a Sasakian structure
with the {\it opposite} orientation. 

\obs All normal $CR$ compact 3-manifolds are covered by circle bundles
over a Riemann surface \cite{gei}, \cite{kamish}, \cite{F}, see also next Section,
and, if this circle bundle is not the Hopf fibration $S^3\rightarrow
S^2$, then all Sasakian structures are, up to a finite quotient,
regular (\cite{F}, see also next Section). If $M$ is covered by $S^3$,
then any Sasakian structure on $M$ is a deformation of a regular one
\cite{F}, see next Section. Therefore the study of these particular
Sasakian structures is essential to the understanding of compact
normal $CR$ 3-manifolds. 
\medskip

As a direct consequence to Proposition \ref{k}, we have, in this case:

\begin{corr}\label{c1}
  Let $M$ be a compact regular Sasakian 3-manifold, thus \linebreak
  $M\stackrel{\pi}{\rightarrow}\Sigma$ is a $S^1$-bundle over 
  a Riemann surface $\Sigma$ Then $-k$ is
  equal to the (Gaussian)   curvature of $\Sigma$.
\end{corr}
The reason for this is that, if $\nabla^\Sigma$ is the Levi-Civita
connection on $\Sigma$, $\nabla$ is the Tanaka-Webster connection on
$M$, and $\tilde X$ denotes the horizontal lift (in $Q$) of a vector
$X\in T\Sigma$, we have:
$$\widetilde{\nb^\Sigma_XY}=\nb_{\tilde X}\tilde Y,\,\forall X,Y\in T\Sigma.$$

We also get, from the above equality, that $\sq f=\sq^\Sigma f$, for a
function $f$ constant on the fibers of $\pi$ (where $\sq^\Sigma$ is
defined by a relation analogous to (\ref{eq:sq})), so the Tanaka
curvature (\ref{eq:fisk}) has a particularly simple expression in this case:

\begin{prop}\label{crpl}
  Let $M$ be a compact regular Sasakian 3-manifold, fibered over a
  Riemann surface $\Sigma$ of positive genus. Then the $CR$ structure of $M$
  is flat iff $\Sigma$ has constant curvature.
\end{prop}
\begin{proof} We have to prove that, if $k$ is the Gauss curvature of
  the Riemann surface $\Sigma$, then it satisfies $\sq k=0, \forall X\in
  T\Sigma$ iff $k$ is
  constant (we omit the indices referring to $\Sigma$, as we only use
  the metric, the Levi-Civita connection and the operator $\sq$ on
  $\Sigma$ in this proof). First we prove the following fact: if
  $f:\Sigma\rightarrow\R$, then 
  \begin{equation}
    \label{eq:sqhol}
    \sq f=0,\forall X\in T\Sigma\Leftrightarrow J(df)^\sharp \mbox{ is
    Killing.}
  \end{equation}
We need to prove that $\nb J(df)^\sharp $ is anti-symmetric, thus it
is enough to check that $g(\nb_XJ(df)^\sharp,X)=0,\forall X\in
T\Sigma$, but this is equal to
{\small $$-X.JX.f+\nb_XJX.f=\h(\sq f+(X.JX-JX.X-\nb_XJX+\nb_{JX}X).f)=\h\sq
f.$$}

If the genus of $\Sigma$ is greater that 1, it admits no non-zero
Killing vector field. If $\Sigma$ is a torus, non-zero vector fields
vanish nowhere, but $J(df)^\sharp$ should vanish in the critical
points of $f$ (e.g. the maximum points). 
%All we have to prove is that the sectional curvature $k$ of $\Sigma$
%cannot satisfy (\ref{eq:sqhol}). In this case,
%``meridians'' (orbits of $df^\sharp$) would be geodesics of equal
%length $l$, and $J(df)^\sharp$ would induce a Jacobi field along each
%such meridian $\gamma:[0,l]\rightarrow \Sigma$. We have thus the equation:
%$$\nb_{\gp} \nb_{\gp} (Jdf^\sharp)=R(\gp,J df^\sharp)\gp=-k Jdf^\sharp.$$
%where $\gp_x=df^\sharp/|df^\sharp|$. If we denote by $\kappa$ the
%length of $ df^\sharp$, then $\kappa=\dot k$ (derivative along
%$\gamma$), and $k$ satisfies the equation
%$$\dddot k=-k\dot k,$$
%thus
%$$\ddot k=-\h\k^2+a,\; a\in\R.$$
%Now, $\ddot k=\dot\kappa$, and $\dot(df)^\sharp=\dot\kappa
%\dot\gamma$, thus it is {\it positive} in $\gamma(0)$ (as
%$(df)^\sharp_{\gamma(0)}=0$ and it increases) and {\it negative} in
%$\gamma(l)$. But $k^2$
\end{proof}
\obs The only possibility to
find a non-constant function $f:\Sigma\rightarrow\R$ such that $\sq
f=0,\forall X$ is that $\Sigma$ is a sphere admitting a isometric
$S^1$ action. Then $f$ is constant on the orbits of this action and
has as only critical points the poles (0-dimensional orbits) of this
action. 

\section[Vaisman metrics and Sasakian structures]{Vaisman metrics on compact complex surfaces
  and Sasakian structures on compact 3-manifolds}\label{s.k3}

%\iVs\iVrg
\begin{defi}
A compact $CR$ 3-manifold $M$ is called {\em primary} if its
fundamental group $\pi_1(M)$ contains no non-trivial finite subgroups.
\end{defi}
This notion arises from the geometry of complex surfaces (see below). 
All circle bundles over a Riemann surface (regular Sasakian manifolds)
are primary; this is no longer the case if we factor them by a finite
group of orientation-preserving bundle automorphisms, unless it acts
either trivially, or with no fixed points, on the basis.   
\smallskip

We recall that a Riemannian product of a Sasakian manifold with a
circle is a Vaisman metric (or locally conformally K{\"a}hler metric
with parallel Lee form) on the resulting complex manifold (which is
usually called a {\it generalized Hopf manifold})
\cite{vais}, \cite{go}, see also \cite{F}. Starting from a compact 3-manifold, we
obtain thus a 
compact complex surface, on which the Vaisman metrics are classified
in \cite{F}:
\begin{prop}
  \label{p1sk} Let $(M,g,T)$ be a Sasakian compact 3-manifold; then
  $M\times S^1$ is one of the following complex surfaces, and $T$ is
  (up to a constant) the following holomorphic vector field:
  \begin{enumerate}
  \item $M\times S^1$ is a properly elliptic surface, admitting two
    holomorphic circle actions: the first one is given by the factor
    $S^1$, and the second one is infinitesimally induced by $T$ (whose
    orbits in $M$ are, therefore, closed). If $M$ is {\rm primary}, then
    the Vaisman metric on $M\times S^1$ is {\rm regular} (i.e. it is
    obtained by a Kaluza-Klein construction on an elliptic fiber
    bundle --- see above);
\item $M\times S^1$ is a Kodaira surface, admitting two
    holomorphic circle actions; all the other conclusions above still
    hold; 
\item $M\times S^1$ is a Hopf surface of class 1, given by the
  contraction $\ul g\in\mathrm{End}(\C^2)$, $\ul g(x,y):=(\alpha
  x,\beta y)$, with $\alpha,\beta\in\R$, $0<\alpha\leq\beta<1$, and
  the Reeb vector field $T$ is induced by the field $i\log\alpha x\d
  x+i\log\beta y\d y$ on $\C^2$. If $M$ is primary, then $M\times S^1$
  is a primary Hopf surface (of class 1) and $M\simeq S^3$ (in
  particular $\pi_1(M)=0$).
  \end{enumerate}
\end{prop}
We recall that a primary properly elliptic surface is a non-flat
elliptic bundle over a Riemann surface of genus $g>1$, a primary
Kodaira surface is a non-flat elliptic bundle over an elliptic curve,
and a primary Hopf surface is a quotient of $\C^2\smallsetminus 0$ by
the infinite cyclic group generated by a contraction $\ul
g$. Non-primary (or secondary) surfaces above considered are finite
quotients of primary ones.
\begin{proof} The almost-complex structure $J^s$ on $M\times S^1$ is defined
  as follows: $J^s|_Q:=J=\nb T$, and $J^s(V):=T$, $V$ being the unitary,
  oriented, generator of $TS^1$. The product metric is, then, given by
  its K{\"a}hler form 
$$\omega:=\h d\eta-\eta\wedge(\eta\circ J^s).$$
It is easy to prove that $J^s$ is integrable and that
\begin{equation}\label{Lee}
d\omega=-2(\eta\circ J^s)\wedge\omega, 
\end{equation}
such that the {\it Lee form} of the Hermitian metric $\omega$, is, by
definition, $\eta\circ J^s$, thus parallel \cite{vais}, \cite{go}, \cite{F}.

The {\it Lee vector field} $V$ (the metric dual of the Lee form) is a
holomorphic vector field \cite{vais}, and it is automatically given
(up to a positive constant) by the complex geometry of the surface
$M\times S^1$ (\cite{F}), which can only be of the three kinds
enumerated in the Proposition. In the first two  cases, it follows
from the descriptions of Sections 3 and 4 of \cite{F}, see also
Theorem 2 of the same paper, that $JV=T$ always has closed orbits. In
the same cases, $V$ and $T$ generate the tangent space of the fibers
of the elliptic fibration, and the Riemannian product situation can
only occur if generic fibers are biholomorphically
$(S^1,\text{can})\times(S^1,l\cdot\text{can})$ (the first factor
corresponds to the orbits of $V$, and the second (of unknown length $l$)
to the orbits of $T$). If $M$ is primary, then $\pi_1(M\times S^1)$
contains no non-trivial finite subgroup, hence $M\times S^1$ is a
principal elliptic bundle over a complex curve of positive genus
(\cite{maehara}, see also
\cite{F}). As any
Vaisman metric is regular on these surfaces, (\cite{F}), Theorem 2, it
immediately follows that the Sasakian metric on $M$ is regular (in the
primary case).
\smallskip

%\ilee\iL\iscp\iskdu\iskdd\isks
If $M\times S^1$ is a Hopf surface, then it is necessarily of {\it
  class 1} (see \cite{F}),
i.e. it is the quotient of $\C^ 2\smallsetminus 0$ by a group $G$,
  generated by a normal finite subgroup $H$ and by a holomorphic
  contraction $\ul g$ of $\C^2$ as in the
Proposition, in general with $\alpha,\beta\in\C$, $0<|\alpha|\leq
|\beta|<1$). We first need the orbits of $V$ to be closed: as
\begin{equation}\label{VV}
V:= \log\alpha x\d x+\log\beta y\d y,\; \text{ see \cite{F},
  Proposition 8},
\end{equation}
we obtain that $\alpha=|\alpha|\varepsilon_1$ and
$\beta=|\beta|\varepsilon_2$, where $\varepsilon_1,\varepsilon_2$ are
both primary $n$-roots of unity for $n\in\N^*$,
i.e. $\varepsilon_j^k=1\Leftrightarrow k=np, \; p\in\Z$, for $j=1,2$.

Then it is clear that $M$ is a finite quotient of $S^3$. In
the primary case, $M\times S^1$ is the primary Hopf surface
$(\C^2\smallsetminus\{0\})/\Z$, where the action of $\Z$ is generated by $\ul
g$. As the generic orbits of $V$ cross the orthogonal hypersurfaces
(tangent to $V^\perp$) $n$ times, it follows that $n=1$ and $M\simeq
S^3$. As a Hopf surface is a finite quotient of a primary one, we get
$n=1$ in general, as claimed.
\end{proof}
\begin{corr}\label{KK}
In the cases 1. and 2. from Proposition \ref{p1sk}, a primary normal
CR 3-manifold is topologically a circle bundle over a Riemann surface
of positive genus. Moreover, the orbits of any CR Reeb vector field
are the fibers of such a fibration, and therefore have all the same
length.
\end{corr}

%As all regular Vaisman metrics are second type deformations of a standard
%one (see \cite{F}), we have that all regular Sasakian metrics are
%{\it second type deformations} of a {\it standard} one; we recall the
%definitions (see also \cite{F}):

We claim that all Sasakian structures on a given 3-manifold are
deformations of a {\it standard} one; in the following definitions (see
also \cite{F}) we currently consider the Sasakian structure $(M,g,T)$
to be given by its $CR$ structure $J$, together with the contact form
$\eta$ and Reeb vector field $T$:

\begin{defi} (i) A {\rm standard} Sasakian structure is, up to finite
  quotient, a regular structure such that the basis has a metric of constant
  curvature; 

(ii) A {\rm first type deformation} of a Sasakian structure $(M,g,T)$
is a Sasa\-kian structure $(M,g',T')$, where $T'=fT+X_f,\, X_f\in Q$ is
another CR Reeb vector field associated to the same $CR$ structure
$(M,Q,J)$. (see also previous Section);

(iii) A {\rm second type deformation} of a Sasakian structure
$(M,\eta,T)$ is defined by the (deformed) contact form $\eta'$,
with the same $CR$ Reeb vector field $T$, and such that:
\begin{equation}
  \label{eq:def2}
  \eta':=\eta+d\sigma\circ J_0,\; Q':=\ker \eta',\; J':={J_0}|_{Q'},
\end{equation}
where $J_0|_Q:=J$ and $J_0(T):=0$ is an extension
$J_0\in\mathrm{End}(TM)$ of $J$; $\sigma$ is a (small enough) function
on $M$, such that $d\sigma(T)=0$ and $\eta'$ is a contact form.

(iv) A {\rm 0-type deformation} of a Sasakian structure consist in
multiplying $T$ with a positive constant, and keeping the $CR$
structure fixed.
\end{defi}
%\iskdz\iskdu\iskdd\icrnns\iSp
After a first type deformation, the normal oriented vector of $Q$
becomes $T'$, and the metric on $Q$ changes as
\begin{equation}
  \label{eq:1def}
  g'|_Q=f^{-1}g|_Q.
\end{equation}
\obs If $(M,g,T)$ is regular, then a second type deformation consist
in a conformal change of the metric on the basis $\Sigma$, and a
subsequent change of connection in the $S^1$-bundle
$M\rightarrow\Sigma$. We can always obtain, by this procedure, a
metric of constant curvature on $\Sigma$, and the corresponding
Sasakian metric on $M$ is determined only by the choice of a
connection of fixed curvature, i.e. of an element of the affine space modeled on
$H^1(\Sigma,\R)$. Standard primary Sasakian structures are then
determined by this latter choice and a complex (conformal) structure
on $\Sigma$; both these choices are unique if $\Sigma\simeq S^2$.

 On any Hopf surface $\C^2\smallsetminus\{0\}/\Z\ul g \ltimes H$ of
 class 1, a Vaisman metric has been constructed in \cite{go}, see also
 \cite{F}, such that it induces a
first type deformation of the round Sasakian structure on $S^3$ (the
unique standard one). We also know \cite{F} that, for a given surface,
 two Vaisman metrics are second type deformations of each other (and
 they induce second type deformations on the corresponding Sasakian
 manifolds). We have then \cite{F}:
\begin{thrm}\label{sasaki}
  Any Sasakian structure on a
  compact 3-manifold $M$ is a deformation of a standard one; a second type
  deformation if $M$ is a Seifert fibration not recovered by $S^3$, a
  composition a first type and of a second type deformation if $M$ is
  a finite quotient of $S^3$. Moreover, a standard Sasakian structure
  on $S^3$ is unique up to a global rescaling (0-type deformation).
\end{thrm}

%\section{$CR$ flat Sasakian structures on compact 3-manifolds}

\section{Normal $CR$ structures on non-spherical Seifert fibrations}\label{s.k4}

Let $M$ be a a Seifert fibration over a 2-dimensional orbifold
$\Sigma$, of non-zero genus. Then $M$ is a finite quotient of a circle
bundle over a Riemann surface of positive genus. If $M$ is primary,
then it is a circle bundle
over a Riemann surface $\Sigma$ which is not a sphere, and we have
seen, in the previous Section, that Sasakian structures on such
manifolds are regular, i.e. given by a
Riemannian metric $g$ on $\Sigma$ and a connection in the $S^1$-bundle
$M\rightarrow \Sigma$, see also \cite{F}. 

Let us consider now $M$ as a normal $CR$ manifold; we ask then if it
admits more than one Sasakian structure associated to it. In other
words: {\it Can a Sasakian structure on $M$ admit $CR$ infinitesimal
automorphisms other than the Killing field $T$?}

The answer is negative; more precisely:

\begin{thrm}\label{t2}
  A Sasakian structure on a non-spherical Seifert fibration does not
  admit non-trivial first type deformations.
  Equivalently, the connected component of the $CR$ automorphism group
of a normal 
  $CR$ structure on $M$ is isomorphic to $S^1$.
\end{thrm}

\begin{proof}
Consider a Sasakian structure on $M$, with the usual notations. (see
Section \ref{s.k2}) We suppose, with no loss of generality, that $M$ is
primary, thus $M\rightarrow
B$ is a $S^1$-bundle over a complex curve $\Sigma$, such that the Reeb
vector field $T$ is tangent to the fibers, which all are of equal
length $l^T$ (and we can suppose $l^T=1$).   

We want to prove that any function $f$ satisfying $\sq f=0,\,\forall
X\in Q$ is constant. The function $\iif f:\Sigma\rightarrow\R$,
defined by 
%\ifc\iskdu
$$\iif f(x):=\int\limits_{M_x}f\eta,$$
where $M_x:=\pi^{-1}(x)$ is a fiber of the projection $M\rightarrow
\Sigma$, satisfies $\sq^\Sigma \iif f=0$, thus, from (\ref{eq:sqhol}),
it is constant (recall that 
$\Sigma$ is of positive genus). If $f>0$ (which we can assume by
adding a constant to the bounded function $f$), then $T'=fT+X_f$ (see
(\ref{eq:xf})) is another $CR$ Reeb vector field on $M$, that has thus
circular
orbits (see Section \ref{s.k3}) of an equal length $l^{T'}$ (as
$(M,Q,J,T')$ is a primary Sasakian manifold, it is regular, see
Section \ref{s.k3}). Suppose $\iif f\equiv 1$. We can also compute 
$$\iif\,'f^{-1}(x'):=\int\limits_{M'_{x'}} f^{-1}\eta',$$
where $M'_{x'}$ is the orbit associated to $x'$ (a point in the orbit
space $\Sigma'$ of $T'$ --- for topological reasons, the projection
$M\longrightarrow\Sigma'$ is still a principal
$S^1$-fibration). $\iif\,'f^{-1}$ is also a constant, and we would like
to relate it to $\iif f$, and to
the length $l^{T'}$ of the orbits of $T'$; the reason for that is the
following Lemma:
\begin{lem}\label{l1sk} % a). Suppose  $T$ and $T'$ have a common orbit
%   $M_x$, such that $f\equiv 1$ on this orbit. 
% Then $f\equiv 1$ everywhere. 
% b). 
a). If two $CR$ Reeb vector fields coincide on a common orbit, they
coincide everywhere.

b). Suppose that $T$ and $T'$ are two arbitrary $CR$ Reeb vector fields,
and that the length of the orbits of each of them (measured in the
corresponding metrics) is equal to 1. If
$\iif f\equiv 1$ and $\iif\,'f^{-1}\equiv 1$, then $T\equiv T'$ everywhere.
% b). Suppose $(\iif\equiv 1\Rightarrow l^f=1)$ for any positive
% function $f$
% satisfying $\sq f=0,\,\forall X\in Q$. Then any such function $f$ is a constant.
\end{lem}
\obs The point {\it a).} is a particular case of {\it b).}; it is
explicitly formulated as it will often be used throughout the proof of
Theorem \ref{t2}.

% The condition $f\equiv 1$ on a common orbit $C$ is equivalent to the
% fact that this orbit has length 1 {\it for both Sasakian metrics}
% Co-responding to $T$ and $T'$~: indeed, as these lengths are $\int
% _C\eta$, resp. $\int_C\eta'$, the second is equal to 
% $$\int_C f^{-1}\eta\ge\left(\int_C\eta\right)^2\left(\int_C
%   f\eta\right)^{-1}=1,$$
% from the H\"older inequality, with equality if and only if $f$ is
% constant along $C$.

\begin{proof} The volume form of the original Sasakian metric is
  $\lambda:=\h\eta\wedge d\eta$, for $\eta=g(T,\cdot)$, and the volume form of
  the deformed metric is $\lambda':=\h\eta'\wedge d\eta'$, where
  $\eta'=f^{-1}\eta$ (see (\ref{eq:1def})). We denote by $v$ the volume
  of $\Sigma$, equal (as 
  the fibers of the $S^1$-bundle $M\rightarrow \Sigma$ are of length 1) to
  the volume of $(M,g,T)$. The volume $v'$ of $(M,g',T')$ is then equal
  to
$$v'=\int\limits_{M}f^{-2}\lambda.$$
On the other hand, from H{\"o}lder's inequality, we have
\begin{equation}
  \label{eq:hold}
  \left(\int\limits_{M} f^{-2}\lambda\right)\left(\int\limits_{M}
  f\lambda\right)^2\geq\left(\int\limits_{M}\lambda\right)^3 = v^3.
\end{equation}
But the integral of $f$
on $M$ is equal to the integral of $\iif f$ on 
$\Sigma$, thus to $v$. We get then
\begin{equation}
  \label{eq:ineg}
  v'\geq v.
\end{equation}
We have thus:
\begin{equation}
  \label{eq:main}
  (\iif f\equiv 1 \mbox{ and }l^T=1)\Rightarrow v'\geq v.
\end{equation}

% Now, let us define the functional $\mbox{\bf\em{I$\,'$}}$ as above, but starting
% from the regular Sasakian structure $(M,g',T')$, and apply it to
% $f^{-1}$. It represents the mean value on the orbits of $T'$ of the
% function $f^{-1}$, which satisfies $\sq' f^{-1}=0,\,\forall X\in Q$
% (and defines the first type deformation of $(M,g',T')$ into
% $(M,g,T)$).

In each of the claims {\it a).} or {\it b).}, both $T$ and $T'$ have orbits
of length 1. On the other hand, we also have
$\iif f\equiv 1\equiv\mbox{\bf\em{I$\,'$}}f^{-1}$, thus we can apply the
implication (\ref{eq:main}) for $f^{-1}$, too (starting from
$(M,g',T')$):
$$(\iif'f^{-1}\equiv 1 \mbox{ and }l^{T'}=1)\Rightarrow v\geq v'.$$
We have thus equality in
(\ref{eq:hold}), which implies that $f$ is constant.

% In case {\it b).}, the condition $l^f=1$ implies the hypothesis of
% (\ref{eq:main}), so we have:
% \begin{equation}
%   \label{eq:11}
%   l^f=1\Rightarrow v'\geq v,
% \end{equation}
% and we can obviously apply it to the Sasakian structure $(M,g',T')$,
% and to the positive function $f^{-1}$. We conclude again that $v'=v$,
% thus $f$ is constant.
\end{proof}

% The strategy to prove Theorem \ref{t2} is to apply part {\it b).} of
% the above Lemma; for this we need to prove that, if $\iif f\equiv 1$,
% then $l^f=1$; we will make use of part {\it a).} of the Lemma to do
% that.
% \smallskip

It is known that the group $G$ of $CR$ automorphisms of $M$ is a Lie group
\cite{kob}, \cite{tk1}. We consider the connected subgroup $G^f$
associated to the Lie subalgebra generated by $T$ and $T'$. This group
acts on $M$ by $CR$ automorphisms, and its orbits are connected.
We consider the decomposition of $M$ in orbits of $G^f$; they are of
three kinds, according to their dimension:
\begin{enumerate}
\item circles $M_{x_i}, \, i\in I\subset B$;
\item immersed surfaces $S\subset M$;
\item open orbits.
\end{enumerate}
% If there is a circular orbit in $M$, then it follows directly from
% part {\it a).} of Lemma \ref{l1sk} that $f\equiv 1$ (in particular
% $l^f=1$, as claimed).

Suppose $S$ is a $G^f$-orbit of dimension 2 in $M$; it contains all
the circles $M_x$ that intersect it, and it is immersed, hence it
projects, via the bundle projection $M\rightarrow \Sigma$, onto an immersed
connected curve $C\subset \Sigma$. There are two cases:
\begin{enumerate}
\item $C$ is an open segment;
\item $C$ is a circle.
\end{enumerate}
In both cases, we define $X^C$ to be a unitary continuous vector field in $TC$,
and $X^S$ to be its horizontal lift to $M$; we have then
$X_f=T'-T=kX^S$, and, as $X_f=\h J(df|_Q)^\sharp$ from (\ref{eq:xf}), we get
\begin{equation}
  \label{eq:xs}
  df(X^S)\equiv 0,
\end{equation}
thus $f$ is constant on the orbits of $X^S$.

% In case 1, these orbits are the fibers of a line bundle
% $S\stackrel{\varphi}{\rightarrow}M_{x_0}$, for a fixed point $x_0\in
% C$, and any orbit $O'$ of $T'$ projects diffeomorphically, via
% $\varphi$, onto $M_{x_0}$, and $T'$ projects onto $fT$. 
% \begin{center}\input{Sasaki/Ssakd.pstex_t}\end{center}
% We see then
% that the length of $O'$ is nothing but
% $$\int\limits_{M_{x_0}}f\eta=1,$$
% thus $l^f=1$ in this case.
% \smallskip

If $C$ is a circle, then it has a finite length $l_C$. Assume that $f$
is not constant on $S$; then $f$ has a regular value $s_0$. Fix
a point $y_0\in M_{x_0}$, such that $f(x_0)=s_0$, and fix some other
points $y_s\in M_{x_0}$ close to it, such that $f(y_s)=s$ are still
regular values of $f$ on $S$. All these
numbers $s$ close to $s_0$ are images of compact immersed curves, 
all diffeomorphic and {\it horizontal} (i.e. tangent
to $X^S$). Consider $C_s$ to be their connected component containing
the points $y_s\in M_{x_0}$; then $C_s$ are precisely the orbits of
$X^S$ starting from $y_s$; they all have thus equal length (which is
an integer multiple of $l_C$, as they all cover $C$). 
%\begin{center}\input{Sasaki/Ssakt.pstex_t}\end{center}

Note that this
length is measured by the Riemannian metric $g$, and that for $g'$ the
lengths need not be equal any longer: Indeed, the curves $C^S$ are
still the same, but the unitary vector fields on them, for $g'$, are
$\sqrt f X^S$ (see Section \ref{s.k3}), and $f$ is precisely non-constant in a
neighbourhood of $y_0$. 

On the other hand, the group $G^f$, the orbit $S$ and the circles
$C_s$ can also be considered starting from the Sasakian metric
$(M,g',T')$, in which case an analog reasoning yields that they have
equal length, in the metric $g'$ (the regular points of $f$ are also
regular points for $f^{-1}$), which leads to a contradiction.

We obtain thus $f\equiv 1$ on $S$. We can suppose that $T'$ is
sufficiently close to $T$ (by replacing $f$, if necessary, with
$\varepsilon f+ 1-\varepsilon$, for $\varepsilon>0$ small), such that
the orbits of $T$ and the ones of $T'$ are homotopic circles in the
torus $S$. We can even suppose that one particular orbit $M'_{x'_0}$
of $T'$ lies 
in a small tubular neighbourhood of an orbit $M_{x_0}$ of $T$. The
projection of this neighbourhood onto $M_{x_0}$ induces then a
diffeomorphism $\phi:M'_{x'_0}\rightarrow M_{x_0}$, such that
$\phi^*\eta=\eta'$ (because $f\equiv 1$ on $S$). Thus 
$$l^{T'}=\int_{M'_{x'_0}}\eta'=\int_{M_{x_0}}\eta=l^T=1,$$
thus we can apply the point {\it b).} in Lemma \ref{l1sk}.
This implies that $f\equiv 1$ (resp. $\varepsilon f+
1-\varepsilon\equiv 1$, which is the same thing).
\smallskip

The same reasoning can be applied if the projection $C$ of $S$ in
$\Sigma$ is a segment of finite length for the metric induced by $T$
(and also for the metric induced by $T'$, if we suppose $T$ and $T'$
to be sufficiently close to each other).
\smallskip

Suppose now $C$ is an immersed open curve in $\Sigma$ of infinite length.
We have 
$$T'=fT+X_f=fT+rX^S,$$
where $r:S\rightarrow\R$ is a function, equal to $-1/2
df(JX^S)$ \ref{eq:xf}. In the following lines, $X$ stands for
$X^S$; it is unitary and contained in $TS\cap Q$:
$$X.JX.f-\nabla_XJX.f-(JX.X.f-\nabla_{JX}X.f)=2T.f=2f',$$
$$X.JX.f-\nabla_XJX.f+(JX.X.f-\nabla_{JX}X.f)=\sq f=0,$$
thus, by summation, we get $X.r-\nabla_XJX.f=f'$, but $\nabla_XJX\perp
JX$ and it still lies in $Q$, hence it is collinear with $X$; but we
know from (\ref{eq:xs}) that $X.f\equiv 0$. We have thus 
$$X.r=f'.$$
We recall that $f$, thus $f'$, too, are constant on the orbits of
$X=X^S$, which are curves of infinite length. If $f'\ne 0$ on such an
orbit, the function $r$ satisfying the equation above is not bounded,
but this is impossible as $M$ is compact.

Thus $f$ is constant on $S$. The rest of the argument used for the
case when $C$ had finite length can also be applied here.
\medskip

We have thus proven that $G^f$ cannot admit any 2-dimensional
orbits. We are going to prove now that the number of 1-dimensional orbits
(which are vertical circles) is finite.

Suppose we have an infinite number of such orbits; then we can extract
a sequence of points $a_n\rightarrow a\in\Sigma, \
n\rightarrow\infty$, such that $M_{a_n},\ n\in\N$, are orbits of
$G^f$. Then so is $M_a$, as the union of all 1-dimensional orbits in
$M$ is closed. We want to prove that $f$ is constant on $M_a$.

The normal bundle of $M_a$ in $M$ is also the
restriction to $M_a$ of the plane bundle $Q$, and will still be
denoted by $Q$.  It is a complex line
bundle, and its metric depends on the $CR$ Reeb vector field in
$\mathfrak{g}^f$ that induces a Sasakian metric on $M$. Consider two
such vector fields $T$ and $T'=fT+X_f$ (where $X_f\equiv 0$ on
$M_a$), and fix two arbitrary different points $x_1,x_2\in M_a$. The
$S^1$ (integral) actions corresponding to $T$, resp. $T'$ on $M$ 
induce two different diffeomorphisms $\psi$,
resp. $\psi':S_1\rightarrow S_2$, where $S_i$ are
contractible surfaces
in $M$, locally defined around $x_i$, tangent to $Q_{x_i}$ at $x_i$
and transverse to the orbits of $G^f$. The choice of such surfaces is
not essential, as we are interested in the differential at $x_i$ of
$\psi$, resp. $\psi'$.

Consider $\Psi:=(\psi')^{-1}\circ\psi$, which is a diffeomorphism of a
neighbourhood $U$ of $x_1$ in $S_1$ into $S_1$. If we denote by $b_n$
the intersections of $M_{a_n}$ with $U$, we obviously have
\begin{equation}
  \label{eq:norm}
  \Psi(b_n)=b_n,\ n\in\N.
\end{equation}
\begin{lem}\label{evid} If $\Psi:U\rightarrow\R^k$ is a
  diffeomorphism from a neighbourhood of $0$ in $\R^k$ into $\R^k$
 that has a sequence, converging to $0$, of fixed points, then the
 differential $(\mathrm{d}\Psi)_0$ has at least one eigenvector
 corresponding to the eigenvalue $1$.
\end{lem}
\begin{proof} By subtracting the inclusion $\mathbf{1}_U$ of $U$ in
 $\R^k$, we get a function $\Psi-\mathbf{1}_U$ which has a sequence,
 converging to $0$, of zeros. Then the kernel of its differential at $0$
 is non-trivial.
\end{proof}
From (\ref{eq:norm}) and the previous Lemma, we conclude that it
exists a non-zero $Y\in Q_{x_1}$ such that
\begin{equation}\label{eq:psi}
(\mathrm{d}\psi)_{x_1}(Y)=(\mathrm{d}\psi')_{x_1}(Y).
\end{equation}
But these differentials are equal to the differentials of the $S^1$
(integral) actions induced by $T$, resp. $T'$, on $M$, and these
actions preserve the complex structure of $Q$. Then $JY$ satisfies
(\ref{eq:psi}) as well, so the two differentials
$(\mathrm{d}\psi)_{x_1}$ and $(\mathrm{d}\psi')_{x_1}$ coincide. On
the other hand, the $S^1$ action induced by $T$ preserves the
corresponding metric $g$ on $Q$, and hence so does
$(\mathrm{d}\psi)_{x_1}$. The same holds for $(\mathrm{d}\psi')_{x_1}$
and $g'=fg$, so we get $f(x_1)=f(x_2)$. As $x_1,x_2$ were arbitrarily
chosen, $f$ is constant on $M_a$, thus everywhere (Lemma \ref{l1sk},
{\it a).}).
\medskip

So the only remaining situation is when $M$ is a union of open orbits
and a finite number of circular orbits of $G^f$; 
as $M$ is connected, there needs to be only one open
orbit $U$ (dense in $M$). We study
now the structure of the Lie algebra $\mathfrak{g}^f$. We will suppose
that $G^f$ acts effectively on $M$.

Every element $V\in\mathfrak{g}^f$ can be written as 
$$V:=f^VT+X_V,$$
where $f^V:M\rightarrow\R$ is a function, and $X_V:=X_{f^V}$. Because 
$\iif f^V$ is a constant, we get a linear homomorphism
$$\iif :\mathfrak{g}^f\rightarrow\R$$
induced by the integral of the functions $f^V$ along the fibers of
$T$. The kernel of this homomorphism is a hyperplane
$H\subset\mathfrak{g}^f$, and it contains all the brackets $[T,V]$,
for any $V\in\mathfrak{g}^f$; indeed, the function $f^{[T,V]}$ is
precisely the derivative along $T$ of $f^V$, denoted by ${f^V}'$,
hence its integral on the orbits of $T$ vanishes.

\begin{lem}
  1. The bracket with $T$,
     $\mathrm{ad}_T\in\mathrm{End}(\mathfrak{g}^f)$, induces an
     automorphism of $H$ (still denoted by $\mathrm{ad}_T$), which is
     $\C$-diagonalizable, and whose eigenvalues are pure imaginary
     (hence non-zero);

\noindent 2. If $\,V$ is the real part of an eigenvector of
     $\mathrm{ad}_T$, then the bracket $[V,\mathrm{ad}_TV]$ is a
     non-zero multiple of $T$. 
\end{lem}
\begin{proof} {\it 1. } Because all orbits of $T$ have length $1$, it
  means that the exponential of $T$, $\exp T\in G^f$, is contained in
  the isotropy subgroup of any point in the open orbit $U$, but the
  intersection of all these isotropy groups is trivial, as $G^f$ acts
  effectively on $M$. In particular, $\mathrm{Ad}_{\exp
  T}=\exp(\mathrm{ad}_T)$ acts trivially on $\mathfrak{g}^f$, so the
  exponential of the endomorphism $\mathrm{ad}_T\in\mathrm{End}(H)$ is
  the identity. It follows that its eigenvalues are imaginary (integer
  multiples of $2\pi i$), and
  that its Jordan decomposition reduces to the diagonal part.

On the other hand, we know from Proposition \ref{crpl} that the only
$CR$ Reeb vector fields commuting with $T$ are multiples of $T$;
therefore, if $[T,V]=0$ for $V\in H$, then $V=0$ and $\mathrm{ad}_T$ is
non-singular, hence all its eigenvalues are non-zero.

\noindent {\bf Remark. } It follows that $\dim H$ is always even, and
that $ \mathrm{ad}_T$ 
is the (commutative) product of a complex structure $\mathcal{J}$ on
$H$ with a diagonal matrix with real eigenvalues.

{\it 2. } If $V$ is the real part of an eigenvector of $\mathrm{ad}_T$,
then $\mathrm{ad}_T^2(V)$ is a multiple of $V$. On the other hand,
$$[T,[V,\mathrm{ad}_TV]]=[\mathrm{ad}_TV,\mathrm{ad}_TV]+
[V,\mathrm{ad}_T^2V]=0,$$  
so $[V,\mathrm{ad}_TV]$ commutes with $T$, hence it is collinear to
it (see Proposition \ref{crpl}).
\end{proof}
Consider the case when the function $f=f^V$ corresponds to the real part
$V$ of an eigenvector of $\mathrm{ad}_T$. Then $G^f$ has dimension 3,
as its Lie algebra is generated by $T,V$, and $\mathrm{ad}_TV$. We
will obtain a contradiction, hence the Theorem will follow.
\smallskip

Denote by $f'$ the function associated to $V'=\mathrm{ad}T(V)$; we have 
\begin{equation}\label{vv}
V=fT+X_f,\ V'=f'T+X_{f'};\ V'=[T,V],\ X_{f'}=[T,X_f].
\end{equation}
We also have 
\begin{equation}\label{sl2r}
[T,V']=-aV, \ a=-4\pi^2l^2,\ l\in\N^*,
\end{equation}
(where $\pm 2\pi il$ are eigenvalues of $\mathrm{ad}_T$) hence $f$,
restricted to any orbit of $T$, is a  
solution of the differential equation
$f''=-af,$
in particular it is a sinusoid function:
\begin{equation}\label{sinus}
f(s)=k_x \sin (2\pi ls),
\end{equation}
for $s$ an arc length
parameter (for the Sasakian metric induced by $T$) on the fiber $M_x$,
and its only critical points are the maximum and the
minimum.
\smallskip

Let us compute, from (\ref{vv}), the bracket $[V,V']$:
\begin{equation}\label{brak}\begin{array}{lcll}
[fT+X_f,f'T+X_{f'}]&\!\!=&\!\!+\,(Ff''-{f'}^2)T&\!\!+\,fX_{f''}-df(X_{f'})T\\
&&\!\!-\,f'X_{f'}+df'(X_f)T&\!\!+\,[X_f,X_{f'}].
\end{array}\end{equation}
As this has to be a constant multiple of $T$, $kT$, it follows that, on a
circular orbit of $G^f$ (where $df=df'=0$), we have
$$ff''-{f'}^2=k,$$
{\it independently on the circular orbit}. But, from (\ref{sinus}),
$ff''-{f'}^2=-4\pi^2l^2k_x^2$, where $l$ is a global constant, and
$k_x$ depends only on the orbit. It follows then that to all circular
orbits $M_x$ of $G^f$ corresponds the same value of $k_x$, the {\it
  amplitude} of the sinusoid $f|_{M_x}$.
\smallskip

The only critical points of $f$ are then its maximums and its
minimums, obtained only on the circular orbits, with the values $\pm
k_x$; indeed, on the open orbit $U$, $V=fT+X_f$ has to be linearly
independent of $T$, thus $df|_Q\ne 0$. 

The function $f$ has the following properties:
\begin{enumerate}
\item it has only a finite set of critical points;
\item any of these (isolated) points is either a maximum or a minimum.
\end{enumerate}
Then, after deforming $f$ if necessary (in order to get a function
with non-singular Hessian at these critical points), we obtain a {\it
  Morse function} $\varphi:M\rightarrow\R$, with a finite set of $2lm$
critical points (where $m$ is the number of circular orbits of $G^f$), 
which are local maximums or local minimums; The
topology of $M$ is thus obtained by glueing $lm$ points to $lm$
2-cells, which implies, as $M$ is connected, that $l=m=1$ and $M$ is
homeomorphic to $S^3$, which contradicts our hypothesis.
% If the $CR$ structure of $M$ is flat, then
% $M\simeq\Gamma\backslash\Hs$ or $M\simeq\Gamma\backslash SL(2,\R)$
% (where $\Gamma$ is a lattice in the corresponding group), and it is
% known that the only infinitesimal $CR$ automorphisms of these flat
% $CR$ manifolds are the 
%  \cite{bsh}). So $U$ is homogeneous,
% non-flat. A theorem of E. Cartan implies then that $\dim G^f=3$
% \cite{cart}.
% We claim that the Lie algebra $\mathfrak{g}^f$ of $G^f$ is generated
% (as a vector space) by $T,T'$ and $[T,T']$, the Lie bracket of $T$ and
% $T'$; indeed, $[T,T']\not\in\R T+\R T'\subset\mathfrak{g}^f$, because
% $\mathfrak{g}^f$ is the Lie algebra generated by $T$ and
% $T'$. In particular, as $M$ is a quotient of $G^f$ by a discrete
% subgroup, $T$ and $T'$ have to be everywhere linearly
% independent. But, in a critical point $x_0$ of $f$ (a maximum point,
% for example), $T'=fT+J(df)^\sharp=fT$, see (\ref{eq:xf}), which is a
% contradiction.
% Thus, for any $CR$ Reeb vector field $T'$, related to $T$ by
% (\ref{eq:xf}), we have $\iif f\equiv 1\Rightarrow l^f=1$. We can apply
% then Lemma \ref{l1sk} to conclude that $f\equiv 1$.
\end{proof}
\begin{corr}
  On a compact, normal $CR$ manifold $M$, the only solutions of the
  equation 
$$ \sq f=0,\;\forall X\in Q$$
are the constants.
\end{corr}
\begin{corr}\label{cor1}
  A compact, normal $CR$ manifold $M$, admits a unique
  Sasa\-kian structure associated to it. For $M$ the total space of a
  circle bundle over a Riemann surface $\Sigma$ of positive genus, an
  isomorphism class of normal $CR$ structures on $M$ is determined
  by an isometry class of Riemannian metrics on $\Sigma$, of unitary
  volume, together with a choice of an element in the affine space
  modeled on $H^1(\Sigma,\R)$.
%For $M$ the total space of a
%  circle bundle over a Riemann surface $B$ of positive genus, the
%  isotopy classes of normal $CR$ structures on $M$ are equivalent
%  to the isotopy classes of Riemannian metrics on $B$.
\end{corr}

%\section{Normal $CR$ structures on $S^3$}

\bigskip
\begin{center}
{\sc {Centre de Math{\'e}matiques\\ UMR 7640 CNRS\\ Ecole Polytechnique
\\91128 Palaiseau cedex\\France}}\\e-mail: {\tt
belgun\@@math.polytechnique.fr}
\end{center}
\bigskip

\end{document}